\documentclass[10pt]{amsart}
\usepackage{amssymb}
\usepackage{amsmath}
\usepackage{subfigure}
\usepackage{amsfonts}
\usepackage{graphicx}
\usepackage[dvips]{epsfig}
\usepackage{geometry}
\usepackage[english]{babel}
\newtheorem{theorem}{Theorem}[section]

\theoremstyle{definition}
\newtheorem{defn}{Definition}

\def\({\left(}

\def\){\right)}

\begin{document}
%%%%%%%%%%%%%%%%%%%%%%%%%%%%%%%%%%%%%
%%%%%%%%%%%%%%%%%%%%%%%%%%%%%%%%%%%%%

\title[Simplicity of the the first eigenvalue   of  $(p,q)$  nonlinear elliptic system ]{Simplicity of the the first eigenvalue   of  $(p,q)$  nonlinear elliptic system }

\author{Farid Bozorgnia}
\address{Department of Mathematics, Instituto Superior T\'{e}cnico, Lisbon.}
\email{bozorg@math.ist.utl.pt}
\thanks{F. Bozorgnia  was supported by the UT Austin-Portugal partnership through
the FCT post-doctoral fellowship
SFRH/BPD/33962/2009 }

\date{\today}
\begin{abstract}
In this short note,  the simplicity of the first eigenvalue of a nonlinear system is shown by an alternative proof; thereby, it states that the first eigenfunctions are unique up to modulo scaling.
\end{abstract}
\maketitle

\section{Introduction}
In  recent years  the eigenvalue problems for  nonlinear elliptic system  involving p-Laplace operator
has  been extensively  studied \cite{ NP, BF}.

 Let $\Omega  \subset  \mathbb{R}^d$ be  a connected, bounded and  open domain  with regular boundary   $\partial\Omega.$  This paper is devoted to show the simplicity of the  first eigenfunction $(u,v) \in W^{1,p}_{0}(\Omega) \times  W^{1,q}_{0}(\Omega)$   of the following eigenvalue system

\begin{equation}\label{1}
 \left\{
\begin{array}{lrl}
 - \Delta_{p} u =\lambda  |u|^{\alpha-1}  |v|^{\beta-1} v  &    \text{in } \Omega, \\
   - \Delta_{q} v =\lambda  |u|^{\alpha-1} |v|^{\beta-1} u  &   \text{in } \Omega,
\end{array}
\right.
\end{equation}
where  $ p, q > 1 $  and $ \alpha, \beta>0 $   are real numbers satisfying
\[
\frac{\alpha}{p}+ \frac{\beta}{q}=1.
\]
Here $\Delta_{p} u =  \text{div}(| \nabla u|^{p-2} \nabla u).$
The  first eigenvalue $ \lambda_{1} (p,q)  $  of system (\ref{1})  is defined  as the least real parameter $ \lambda $  for
which both equations of (\ref{1}) have a nontrivial solution $(u,v)$   in
  $ W^{1,p}_{0}(\Omega) \times W^{1,q}(\Omega),$ where  with $ u,  v \neq 0.$

The coupled system    (\ref{1})  has several applications. For instance, in the case where $ p > 2,$  problem (\ref{1})  appears in the study of
non-Newtonian fluids, pseudoplastics for $1 < p <2,$  and in reaction-diffusion problems,
flows through porous media, nonlinear elasticity, and glaciology for $p = \frac{4}{3}.$ We refer the interested readers to     \cite{Diam} for more details.

In \cite{AMZ,KMO} it is shown that  the first eigenvalue
is  simple and corresponding eigenfunction is  non negative. Moreover, the stability  of the first eigenvalue with respect to $(p,q) $  is established. The proof given in this work is  simpler  compare with mentioned works.

\section{ The first eigenvalue of non linear elliptic system  }

Consider   the following   system

\begin{equation}\label{3}
 \left \{
\begin{array}{lrl}
 - \Delta_{p} u =\lambda  |u|^{\alpha-1}  |v|^{\beta-1} v  &    \text{in } \Omega, \\
   - \Delta_{q} v =\lambda  |u|^{\alpha-1} |v|^{\beta-1} u  &   \text{in } \Omega,\\
   u=v=0                &   \text{on  } \, \partial\Omega,
\end{array}
\right.
\end{equation}
where  $ p, q > 1$  and $ \alpha, \beta>0 $   are real numbers satisfying
\[
\frac{\alpha}{p}+ \frac{\beta}{q}=1.
\]
\begin{defn}
The  first eigenvalue  $\lambda_{1}(p,q)$  of (\ref{3})   is defined  as the least real parameter $ \lambda$ for
which both equations of (\ref{3})   have a nontrivial solution $(u,v)$  in the product Sobolev space $ W^{1,p}_{0} \times W^{1,q}_{0} $  with $ u\neq 0$  and $ v\neq  0. $
\end{defn}
  Here by a solution to (\ref{3}) we mean a pair $(u,v)$ such that
  \[
\int_{\Omega}|\nabla u|^{p-2} \nabla u \cdot \nabla\phi \, dx  +  \int_{\Omega}|\nabla v|^{q-2} \nabla u \cdot \nabla\psi  \, dx =\lambda( \int_{\Omega}| u |^{p-2}  u \, \phi \, dx + \int_{\Omega}| v |^{p-2}  v \, \psi \, dx),
\]
\[  \quad \forall \, (\phi, \psi )  \in
 W^{1,p}_{0}(\Omega) \times  W^{1,q}_{0}(\Omega).
\]
The principal eigenvalue  $\lambda_{1}(p, q)$   can be variationally characterized by  minimizing the functional
\[
I(u, v) = \frac{\alpha}{p}\int_{\Omega} |\nabla u(x)|^{p} dx  + \frac{\beta}{q}\int_{\Omega} |\nabla v(x)|^{q} dx,
\]
over the set
\[
C={\{(u, v) \in  W^{1,p}_{0} \times  W^{1,q}_{0}: \int_{\Omega} |u(x)|^{\alpha-1}|v (x)|^{\beta-1} u(x) v(x) dx =1 }\},
\]
by definition
\begin{equation}\label{4}
\lambda_{1}(p,q) =\inf{\{ I(u,v), (u,v) \in C}\}.
\end{equation}
The pair $(u,v) $  is called an eigenvector.
Note that that the solutions  $(u,v) $  of  (\ref{3})   correspond to the critical points of the energy functional $I (u,v).$
%For different choices of $p, q, \alpha, \beta $ problem (\ref{3}) has many  applications for instance in  reaction-diffusion problems, flows through porous media, nonlinear elasticity, and  non-Newtonian fluids.
%%
%\begin{remark}
%The first eigenvalue $ \lambda_{1}(p,q) $ defined by (\ref{4}) is simple an corresponding eigenfunction  functions  $(u,v) $ is nonnegative.
%In \cite{}, the existence of a principal eigenvalue and simplicity and isolation of the principal eigenvalue also positivity of associated eigenfunction are shown. Here we  give an alternative proof for simplicity of first eigenvalue of system in (\ref{3}).
%
%\end{remark}

\begin{theorem} Let $ \lambda_{1}(p,q) $  be defined by (\ref{4}), then   $ \lambda_{1}(p,q)$ is  simple.

\begin{proof}
Let  $(u,v) $ and  $ (\phi,\psi ) $   be two    normalized  vectors of   the first eigenfunction   associated with  $ \lambda_{1}(p,q) $.  We show
that there exist real numbers $ k_1, k_2 $ such that
$ u=k_1 \phi $ and  $ v=k_2 \psi $.
Our proof is based on proof given by Belloni and Kahwol (see \cite{BK}  for scaler   case).  Note that  the function defined below are admissible functions for problem (\ref{4}):
\[
w_1=(\frac{u^{p}+\phi^{p}}{2})^{\frac{1}{p}} \quad \text{and} \quad   w_2=(\frac{v^{q}+\psi^{q}}{2})^{\frac{1}{q}}.
\]
We have
\[
|\nabla w_1|^{p}= (\frac{u^{p}+\phi^{p}}{2}) |\frac{ u^{p} \nabla \log \, u + \phi^{p}  \nabla \log\, \phi}{u^{p}+\phi^{p}}|^{p},
\]
\[
|\nabla w_2|^{q}= (\frac{v^{q}+\psi^{q}}{2}) |\frac{ v^{q} \nabla \log \, v + \psi^{q}  \nabla \log\, \psi}{v^{q}+\psi^{q}}|^{q}.
\]
Now by Jensen's inequality  for convex function $  \theta(\cdot)= |\cdot |^{p}$  with
\[
\theta(\frac{\Sigma a_i x_i}{\Sigma a_i}) \le  \frac{\Sigma a_i \theta(x_i)}{\Sigma a_i},
\]
choose

\begin{equation*}
 \left \{
\begin{array}{ll}
a_1=  \frac{u^{p}}{u^{p}+\phi^{p}}, \quad  a_2= \frac{\phi^{p}}{u^{p}+\phi^{p}},\\
x_1= \nabla \log\, u, \quad x_2= \nabla \log\, \phi.
\end{array}
\right.
\end{equation*}
Thus the following inequalities  hold
\[
|\nabla w_1|^{p} \le \frac{1}{2} |\nabla  u|^{p} + \frac{1}{2}|\nabla  \phi|^{p},
\]
\[
|\nabla w_2|^{q} \le \frac{1}{2} |\nabla  v|^{q} + \frac{1}{2}|\nabla  \psi|^{q}.
\]
The inequalities above are strict at points where  $ \nabla \log \,  u\neq \nabla \log \, \phi $ and  $ \nabla \log \,  v\neq \nabla \log \, \psi.$ Now we have
\begin{equation}\label{f23}
\lambda_{1}\le \frac{\frac{\alpha}{p}\int_{\Omega} |\nabla w_1|^{p} dx  + \frac{\beta}{q}\int_{\Omega}
 |\nabla w_2|^{q} dx}{\int_{\Omega}w_{1}^{\alpha} w_{2}^{\beta} \, dx}.
\end{equation}
Note that by concavity we have
\[
w_{1}^{\alpha} w_{2}^{\beta} =  (\frac{u^{p}+\phi^{p}}{2})^{\frac{\alpha}{p}} \cdot  (\frac{v^{q}+\psi^{q}}{2})^{\frac{\beta}{q}}\ge (\frac{u^{\alpha}+\phi^{\alpha}}{2}) \cdot  (\frac{v^{\beta}+\psi^{\beta}}{2}).
\]
It is easy to see that $ u, v, \phi, \psi $ can be chosen such that
\[
\int_{\Omega} u ^{\alpha} \, v ^{\beta}\,dx =\int_{\Omega} u^{\alpha} \, \psi^{\beta}\,dx=\int_{\Omega} \phi^{\alpha} \, v ^{\beta}\,dx=\int_{\Omega} \phi^{\alpha} \, \psi^{\beta}\,dx=1.
\]
Then  the inequality in (\ref{f23}) reads as
\begin{equation}\label{f24}
\lambda_{1}\le \frac{\frac{\alpha}{2p} (|\nabla  u|^{p} + |\nabla  \phi|^{p} )  + \frac{\beta}{2q}( |\nabla  v|^{q} + |\nabla  \psi|^{q}) } {\frac{1}{4} \int_{\Omega}  u ^{\alpha} \, v ^{\beta} +  u^{\alpha} \, \psi^{\beta} + \phi^{\alpha} \, v ^{\beta}+ \phi^{\alpha} \, \psi^{\beta}  \, dx}.
\end{equation}
If  $ \nabla \log \,  u\neq \nabla \log \, \phi $ and  $ \nabla \log \,  v\neq \nabla \log \, \psi $ in a set
 of positive measure, then we would have strict inequality above, which is contradiction.  This shows that $u$ and $v$ also $\phi$ and  $\psi $ are constant multiplies of each other.
\end{proof}
\end{theorem}
%\section{Algorithm  for a class of quasi linear Elliptic systems}
%
%The algorithm can be applied for the following system
%\begin{equation}\label{f6}
% \left\{
%\begin{array}{lrl}
% - \Delta_{p} u =\lambda     F_{u}(x,u,v) &    \text{in } \Omega, \\
%   - \Delta_{q} v =\lambda   F_{v}(x,u,v)   &   \text{in } \Omega,\\
%   u=v=0     &   \text{on }   \partial\Omega,\\
%\end{array}
%\right.
%\end{equation}
%where  $  1 < p, q < \infty. $
% The   function $  F (x,u,v) $   satisfy the following
% \begin{itemize}
%\item   $ F(x,  0, 0) =  F_{u}(x, 0,0) = F_{v}(x, 0,0) =0$
%\item   $| F_{u}(x, u,v)| \le C(1+ |t|^{p-1} + |s|^{q\frac{p-1}{p}} ) $
%\item   $ | F_{v}(x, u,v)| \le C(1+ |t|^{q-1} + |s|^{p \frac{q-1}{q}} ) $
%\end{itemize}
%Under the growth assumption on $F$ the weak solutions of  (\ref{f6}) are critical points of the following functional
%Consider
%\[
%I(u, v) = \frac{\alpha}{p}\int_{\Omega} |\nabla u(x)|^{p} dx  + \frac{\beta}{q}\int_{\Omega} |\nabla v(x)|^{q} dx,
%\]
%over the set
%\[
%C={\{(u, v) \in  W^{1,p}_{0} \times  W^{1,q}_{0}: \int_{\Omega} F(x,u(x), v(x)) \,  dx =1 }\}.
%\]
%
% Under this assumptions, the critical points of  $I(u,v)$ is the following elliptic system
%\begin{equation}\label{f7}
% \left\{
%\begin{array}{lrl}
% - \Delta_{p} u =\lambda     F_{u}(x,u,v) &    \text{in } \Omega, \\
%   - \Delta_{q} v =\lambda   F_{v}(x,u,v)   &   \text{in } \Omega,
%\end{array}
%\right.
%\end{equation}
%
% We remark that (2) implies that F satises the following growth condition:
%  $ | F(x,  t,  s)| \le C(1 + |t|^p + |s|^q)  $

%\begin{thebibliography}{widest-label}

{}
%\end{large}
 \end{document}